\documentclass[11pt]{article}

\usepackage[T1]{fontenc}
\usepackage[utf8]{inputenc}
\usepackage{lmodern}
\usepackage{amsmath,amssymb,amsthm,mathtools}
\usepackage{geometry}
\usepackage{microtype}
\usepackage{enumitem}
\usepackage{booktabs}
\usepackage{xcolor}
\usepackage[hidelinks]{hyperref}

\allowdisplaybreaks

\hypersetup{
  pdftitle={Connected graphs having maximum signless Laplacian spectral radius},
  pdfsubject={Signless Laplacian spectral extremal graphs},
  pdfkeywords={signless Laplacian spectral radius, threshold graph, Ferrers diagram, coronal, quasi-star, quasi-complete}
}

\newtheorem{theorem}{Theorem}[section]
\newtheorem{lemma}[theorem]{Lemma}
\newtheorem{proposition}[theorem]{Proposition}
\newtheorem{corollary}[theorem]{Corollary}

\theoremstyle{definition}
\newtheorem{definition}[theorem]{Definition}
\theoremstyle{remark}
\newtheorem{remark}[theorem]{Remark}

\newcommand{\Cnm}{\mathbf C_{n,m}}
\newcommand{\one}{\mathbf 1}
\newcommand{\join}{\mathbin{\vee}}
\newcommand{\disj}{\mathbin{\sqcup}}

\definecolor{revisionblue}{RGB}{20,82,145}

\title{Connected graphs having maximum signless Laplacian spectral radius}
\author{Di-Jian  Wang, Hong-Yan Ge, Shi-Cai Gong\thanks{Corresponding author.   E-mail addresses: dijianwang0121@163.com(D.J. Wang), 398643629@qq.com (H.Y. Ge), and
scgong@zafu.edu.cn(S.C. Gong), ybjmath@163.com(B.J. Yuan).} ~and Bo-Jun Yuan
\\{\small \it  School of Science, Zhejiang University of Science and Technology, }\\{\small \it
Hangzhou, 310023, P. R. China}}

\date{}

\begin{document}
\maketitle

\begin{abstract}
Characterizing extremal graphs is a central problem in spectral graph
theory. Let $\rho(G)$ denote the signless Laplacian spectral radius
of a graph $G$, and let $\mathcal C_{n,m}$ be the family of connected
graphs of order $n$ and size $m$. A graph in $\mathcal C_{n,m}$ is
called $Q$-maximal if it attains the maximum signless Laplacian
spectral radius over this family. Motivated by a question posed by
Chang and Tam in 2011 concerning quasi-star and quasi-complete
extremizers, we investigate the corresponding connected extremal
problem using two complementary tools.

First, we establish an identity based on the Perron equations of a threshold graph  and then
 obtain the scalar identity equation
\[
  S_G(\lambda)=T_{a,b}(\lambda),
  \qquad \lambda=\rho(G),
\]
where $a$ and $b$ are the total numbers of zero-vertices and
one-vertices, respectively. This identity yields a one-point
comparison principle: if a  connected threshold graph $H$ with the
same block totals satisfies
$S_H(z)>T_{a,b}(z)$ for some $z\ge n$, then $\rho(H)>z$.
Consequently, spectral improvement can be established without
estimating individual Perron coordinates. Second, we adapt coronal
methods to threshold graphs, providing an additional
approach to spectral comparison.

Our main structural result shows that every $Q$-maximal threshold
graph has character at most three; equivalently, its vertex-sequence
contains at most three zero--one block pairs. We further prove that
the quasi-star graph attains the maximum signless Laplacian spectral
radius over $\mathcal C_{n,m}$ whenever $
  n\ge5$ and $
  n-1\le m<\frac{n^2+10n-8}{8}.$
Thus quasi-star maximality holds throughout this explicit sparse
range. For certain dense parameter pairs, however, we construct
connected threshold graphs whose signless Laplacian spectral radii
strictly exceed those of both the quasi-star and the connected
quasi-complete graph. Hence the corresponding two-candidate
assertion for connected graphs fails in general.
\end{abstract}



\section{Introduction}\label{sec:introduction}

For a finite simple undirected graph $G$, the
\emph{signless Laplacian matrix}, or \emph{$Q$-matrix}, is $ Q(G)=D(G)+A(G),
$ where $A(G)$ is the adjacency matrix of $G$ and $D(G)$ is the
diagonal matrix of its vertex degrees. The multiset of eigenvalues
of $Q(G)$, counted with multiplicity, is called the
\emph{$Q$-spectrum} of $G$. Two graphs are called
\emph{$Q$-cospectral} if their $Q$-spectra coincide; if they are
also nonisomorphic, they form a nonisomorphic $Q$-cospectral pair.
The largest eigenvalue of $Q(G)$, denoted by $\rho(G)$, is called
the \emph{signless Laplacian spectral radius} of $G$.

Let $\mathcal G_{n,m}$ denote the family of all graphs of order
$n$ and size $m$, and let $\mathcal C_{n,m}$ be its subfamily
consisting of connected graphs. A graph $G\in\mathcal C_{n,m}$
is called \emph{$Q$-maximal} if
\[
  \rho(G)=\max\{\rho(H):H\in\mathcal C_{n,m}\}.
\]

The problem of maximizing the spectral radius for prescribed
order and size has a long history in spectral graph theory.
For the adjacency setting, a foundational contribution is
Schwarz's rearrangement principle for nonnegative
matrices~\cite{Schwarz1964}. Brualdi and Hoffman subsequently
studied spectral-radius maximization for $(0,1)$-matrices and
conjectured that, among simple graphs with a prescribed number
of edges, the maximum adjacency spectral radius is attained by
the ordinary, possibly disconnected, quasi-complete
graph~\cite{BrualdiHoffman1985}. Substantial partial results were
obtained by Friedland and Stanley, and Rowlinson established the
conjecture using a double-eigenvector
argument~\cite{Friedland1985,Stanley1987,Rowlinson1988}.
Consequently, the maximum adjacency spectral radius over
$\mathcal G_{n,m}$ is attained by the ordinary quasi-complete
graph. More recently, Cheng and Weng developed a systematic
approach based on modified quotient matrices and applied it to
characterize spectral-radius-maximizing $(0,1)$-matrices with a
prescribed number of ones~\cite{ChengWeng2025}. These developments
illustrate the effectiveness of rearrangement arguments,
Ferrers-type structures, and quotient-matrix reductions in
spectral optimization.

The connected adjacency problem is considerably more delicate.
In 1986, Brualdi and Solheid investigated the graphs of maximum
adjacency spectral radius in $\mathcal C_{n,m}$ and determined
the extremal graphs in several sparse
ranges~\cite{BrualdiSolheid1986}. Further cases were resolved by
Cvetkovi\'c and Rowlinson and by
Bell~\cite{Bell1991,CvetkovicRowlinson1988}. Nevertheless, a
complete solution for arbitrary $n$ and $m$ remains unavailable.
An important structural distinction is that isolated vertices
can be added to a dense extremal component in the unrestricted
problem, whereas connectedness requires every vertex to be
incorporated into a single component using the same edge budget.
This distinction motivates the consideration of both
clique-concentrated and star-concentrated constructions.

The corresponding signless Laplacian problem has also attracted
considerable attention. Unlike the adjacency spectral radius,
the $Q$-spectral radius reflects both the adjacency pattern and
the redistribution of vertex degrees under an edge
transformation. Early progress concerned graphs of small
cyclomatic number. In particular, Fan, Tam, and Zhou determined
the extremal bicyclic graph~\cite{FanTamZhou2008}.
Chang and Tam subsequently proved that, for
$3\le \ell\le n-3$, the graph obtained from $K_{1,n-1}$ by joining
one leaf to $\ell+1$ other leaves is the unique graph in
$\mathcal C_{n,n+\ell}$ having maximum signless Laplacian
spectral radius~\cite{ChangTam2011}. Zhai, Lin, and Zhao extended
this conclusion, for $4\le \ell\le n-3$, from connected graphs
to all graphs in
$\mathcal G_{n,n+\ell}$~\cite{ZhaiLinZhao2022}.
More recently, Zhang et al.~\cite{ZhangJinWangYangZhang2025}
studied the corresponding extremal problem for $
  A_\alpha(G)=\alpha D(G)+(1-\alpha)A(G),
   0\le\alpha<1,
$
and characterized the extremal graphs in an explicit linear-edge
range: for a positive integer $r$, they considered
$n>30r$ and  $
  n-1\le m\le rn-\frac{r(r+1)}2.$
Since $A_{1/2}(G)=Q(G)/2$, their results also yield further
conclusions for $Q$-maximal graphs.

Before stating the problem, we introduce the two principal
candidate graphs and two related families of threshold graphs.
Throughout, assume that $n\ge2$ and
$n-1\le m\le\binom n2$. We use $\join$ and $\disj$ to denote
the join and disjoint union of graphs, respectively, and
$K_0$ denotes the empty graph on zero vertices.

At the complete-graph endpoint, set $
  C_{n,\binom n2}=S_{n,\binom n2}=K_n.$
For the remaining cases, namely
$n-1\le m<\binom n2$, write uniquely
\[
  m-(n-1)=\binom a2+b,
  \qquad 0\le b<a\le n-2.
\]
The \emph{connected quasi-complete graph} is defined by
$C_{n,m}=K_1\join
  [
    (K_b\join(K_{a-b}\disj K_1))
    \disj(n-a-2)K_1
  ].
$
If $b\le n-1-a$, we  define 
$
  C_{n,m}^{\triangle}
  =K_1\join
  [
    (K_1\join(K_{a-1}\disj bK_1))
    \disj(n-1-a-b)K_1
  ].
$
Obviously,  $C_{n,m}, C_{n,m}^{\triangle}\in \mathcal C_{n,m}$.

Our convention for the quasi-complete graph differs from the
usual one in the literature; see, for
example,~\cite{AhlswedeKatona1978}. If
$\widetilde C_{n,m}$ denotes the ordinary quasi-complete graph
of order $n$ and size $m$, then
$
  C_{n,m}=K_1\join\widetilde C_{n-1,m-n+1}.
$
Thus $C_{n,m}$ is obtained by adding a dominating vertex to the
ordinary quasi-complete graph of order $n-1$ and size $m-n+1$.
Unless otherwise specified, the term \emph{quasi-complete graph}
will refer to this connected version throughout the paper.

For $n-1\le m<\binom n2$, there are unique integers $j,j'$
satisfying
\[
  m=\binom n2-\binom{j+1}{2}+j',
  \qquad
  2\le j\le n-1,\quad 1\le j'\le j.
\]
The \emph{quasi-star graph} is
$ S_{n,m}=K_{n-j-1}\join
    \bigl(K_{1,j'}\disj(j-j')K_1\bigr).$
When $j'=3$ and $3\le j\le n-2$, we define 
$S^\triangle_{n,m}
  =K_{n-j-1}\join
    \bigl(K_3\disj(j-2)K_1\bigr).$
This construction replaces the claw $K_{1,3}$ by
$K_3\disj K_1$, preserving both the order and the size.
The condition $j\le n-2$ ensures that the joining clique is
nonempty and hence that $S^\triangle_{n,m}$ is connected.

Motivated by the question posed by Chang and
Tam~\cite[p.~2458]{ChangTam2011}, we consider the following
connected analogue.

\medskip
\noindent\textbf{Problem 1.}
For every pair of integers $n,m$ satisfying
$n-1\le m\le\binom n2$, is it true that
\[
  \max_{G\in\mathcal C_{n,m}}\rho(G)
  =\max\{\rho(S_{n,m}),\rho(C_{n,m})\}?
\]

We approach this problem through the structural properties of
$Q$-maximal graphs.
A graph is called a \emph{threshold graph} if it contains no
induced subgraph isomorphic to $2K_2$, $P_4$, or $C_4$.
The following result is fundamental.

\begin{lemma}[\cite{FanTamZhou2008}]\label{lemma1}
For every pair of integers $n,m$ satisfying
$n-1\le m\le\binom n2$, every $Q$-maximal graph in
$\mathcal C_{n,m}$ is a threshold graph.
\end{lemma}

A connected threshold graph can be represented by a
\emph{vertex-sequence}
$
  (a_1,\ldots,a_k;\,b_1,\ldots,b_k),$
where $a_i$ and $b_i$ are the sizes of its successive zero-blocks
and one-blocks, respectively. The number $k$ of zero--one block
pairs is called the \emph{character} of the graph.
The representation and its associated Ferrers diagram are
described in detail in Section~\ref{sec:prelim}.

Our first main result bounds the structural complexity of
$Q$-maximal graphs.

\begin{theorem}\label{thm:main1}
Every $Q$-maximal graph in $\mathcal C_{n,m}$ is a threshold
graph with character at most three.
\end{theorem}

Combining this structural reduction with further necessary
conditions, we establish quasi-star maximality in the
following range of edge counts.

\begin{theorem}\label{thm:main2}
Let $n\ge5$ and $n-1\le m<\frac{n^2+10n-8}{8}$. Then
 $S_{n,m}$ is $Q$-maximal.
Moreover, if $j'=3$, the  graph $S^\triangle_{n,m}$ defined
as above is also $Q$-maximal,
and $
  \rho(S^\triangle_{n,m})=\rho(S_{n,m}).$
\end{theorem}

\begin{remark}
  The bound in Theorem~\ref{thm:main2} is sharp in the sense that
the strict inequality cannot be replaced by a non-strict one.
Indeed,
\[
  \rho(S_{6,11})<\frac{33}{4}
  <\rho(C_{6,11})=\frac{9+\sqrt{57}}2,
  \qquad
  11=\frac{6^2+10\cdot6-8}{8}.
\]
\end{remark}
\begin{remark}
Problem~1 has a negative answer in general. For example, let
$
  C^\triangle_{9,25}
   =K_2\join(K_5\disj2K_1).$
Then
$
  \max\{\rho(C_{9,25}),\rho(S_{9,25})\}
  <\frac{90}{7}<\rho(C^\triangle_{9,25}).
$
Thus neither $C_{9,25}$ nor $S_{9,25}$ is $Q$-maximal.
Unlike the companion $S^\triangle_{n,m}$, which shares the
$Q$-spectrum of the quasi-star, this example yields a strict
spectral improvement over both standard candidates.
\end{remark}

The remainder of the paper is organized as follows.
Section 2 introduces the notation and develops the scalar comparison criterion and the \(Q\)-coronal identities. Section 3 states the structural restrictions and spectral estimates used in the proofs of the main results. Section 4 proves Theorems \ref{thm:main1} and \ref{thm:main2}. The proofs of the auxiliary results stated in Section 3 are collected in Section 5.

\section{Preliminaries}\label{sec:prelim}

All graphs considered in this paper are finite, simple, and
undirected. For a graph $G$, its order and size are
$|V(G)|$ and $|E(G)|$, respectively. For a vertex $v\in V(G)$,
let $ N_G(v)$ denote its \emph{open neighborhood}, and let
$d_G(v)=|N_G(v)|$ be its degree.  Its \emph{closed neighborhood}
is $N_G[v]=N_G(v)\cup\{v\}$.
The vertex $u$ is called a {\it dominating vertex} of the graph 
$G$ if $d_G(u)=|V(G)|-1$ and is called an {\it isolated vertex} if $d_G(u)=0.$ 
 When the underlying graph is clear,
we omit the subscript $G$.
For $U\subseteq V(G)$, write $G[U]$ for the subgraph induced
by $U$ and $G-U$ for $G[V(G)\setminus U]$.

Two distinct vertices $u,v\in V(G)$ are called \emph{twins} if
$
  N_G(u)\setminus\{v\}=N_G(v)\setminus\{u\}.
$
They are called \emph{false twins} if they are nonadjacent,
in which case $N_G(u)=N_G(v)$, and \emph{true twins} if they
are adjacent, in which case $N_G[u]=N_G[v]$.

The complement of $G$ is denoted by $\overline G$.
We write $K_r$ for the complete graph of order $r$ and
$K_{s,t}$ for the complete bipartite graph with part sizes
$s$ and $t$. The notation $rK_1$ denotes the edgeless graph
on $r$ vertices.

For matrices, $I_r$ denotes the identity matrix of order $r$,
and $J_{r,s}$ and $O_{r,s}$ denote the $r\times s$ all-ones
and zero matrices, respectively. We abbreviate
$J_{r,r}$ and $O_{r,r}$ as $J_r$ and $O_r$.
The all-ones column vector of length $r$ is denoted by
$\mathbf 1_r$.

For a connected graph $G$, a positive eigenvector of $Q(G)$
corresponding to $\rho(G)$ is called a \emph{Perron vector}.

A threshold graph can be constructed from a single vertex by
successively adding a vertex that is either isolated  or
 adjacent
to every existing vertex. These two operations are
encoded by $0$ and $1$, respectively, and the initial vertex
is encoded by $0$.

For a connected threshold graph of order at least two, the
creation word ends with $1$. Grouping consecutive identical
symbols into maximal runs gives its \emph{vertex-sequence}
\begin{equation}\label{eq:vertex-sequence}
  G=(a_1,\ldots,a_k;\,b_1,\ldots,b_k),
  \qquad a_i,b_i\in\mathbb Z_{>0}.
\end{equation}
The corresponding  binary {\it creation word} is
$
  0^{a_1}1^{b_1}\cdots0^{a_k}1^{b_k},
$
where the exponents indicate the lengths of the respective
runs. We use {\it creation-word} notation only when describing
extensions obtained by appending individual threshold symbols;
otherwise, we use vertex-sequence notation.

Recall that the number $k$ of zero--one block pairs is called
the \emph{character} of $G$. For convenience in transformation
arguments, we occasionally allow blocks of size zero in
temporary block descriptions. Before computing the character,
we delete empty blocks, restore the convention that the initial
vertex is encoded by $0$ if necessary, and merge consecutive
runs of the same symbol. Thus the character is always computed
from the normalized representation with positive block sizes,
rather than from the temporary list of parameters.

Let $A_i$ and $B_i$ denote the vertex sets corresponding to
$0^{a_i}$ and $1^{b_i}$, respectively. Then
$\bigcup_{i=1}^k A_i$ is an independent set,
$\bigcup_{i=1}^k B_i$ induces a clique, and
$
  N_G(u)=\bigcup_{j=i}^k B_j$ 
  for every $u\in A_i.$
Consequently, any two distinct vertices in the same zero-block
are false twins, while any two distinct vertices in the same
one-block are true twins.

For the vertex-sequence in \eqref{eq:vertex-sequence}, set
\begin{equation*}
\begin{aligned}
  a&=\sum_{i=1}^k a_i,
  &\qquad b&=\sum_{i=1}^k b_i,
  &\qquad s_i&=\sum_{j=1}^i a_j,
  & \qquad t_i&=\sum_{j=i}^k b_j,
  \qquad(1\le i\le k),
\end{aligned}
\end{equation*}
with $s_0=0$. Thus $a$ and $b$ are the total numbers of
zero-vertices and one-vertices, respectively.
By the adjacency relations described above,
$$
  d_G(u)=t_i \quad(u\in A_i)~~and ~~
  d_G(v)=b+s_i-1 \quad(v\in B_i).
$$
Since the one-vertices induce a clique and the zero-vertices
form an independent set, we obtain
\begin{equation*}
\begin{aligned}
  n&=a+b,~~
  m=\binom b2+\sum_{i=1}^k b_i s_i
    =\binom b2+\sum_{i=1}^k a_i t_i.
\end{aligned}
\end{equation*}
The two sums count the edges between the zero-vertices and
one-vertices from the two sides of the partition.

Using the Ferrers-diagram viewpoint employed by Luo, Gong,
and Tian~\cite{LuoGongTian2024}, order the vertex blocks as
\begin{equation}\label{eq:ferrers-order}
  (B_k,B_{k-1},\ldots,B_1,A_1,A_2,\ldots,A_k).
\end{equation}
Let $\mathcal F(G)$ denote the $b\times a$ submatrix of $A(G)$
whose rows correspond to the one-vertices and whose columns
correspond to the zero-vertices, in the order specified in
\eqref{eq:ferrers-order}. Its block with row set $B_r$ and
column set $A_j$ is
\begin{equation*}
  \mathcal F(G)[B_r,A_j]
  =
  \begin{cases}
    J_{b_r,a_j},& j\le r,\\
    O_{b_r,a_j},& j>r,
  \end{cases}
  \qquad 1\le r,j\le k.
\end{equation*}
In particular, every row corresponding to a vertex of $B_r$
contains exactly $s_r$ ones, and every column corresponding
to a vertex of $A_j$ contains exactly $t_j$ ones.
With the one-vertices and zero-vertices grouped together,
the adjacency matrix takes the form
\begin{equation}\label{eq:adj-Ferrers}
  A(G)=
  \begin{pmatrix}
    J_b-I_b & \mathcal F(G)\\
    \mathcal F(G)^{\mathsf T} & O_a
  \end{pmatrix}.
\end{equation}
In the chosen order, the ones in each row of $\mathcal F(G)$
are left-justified, and the row lengths are weakly decreasing.
Thus $\mathcal F(G)$ is a rectangular Ferrers matrix, expressing
the nested neighborhoods between the one-vertices and
zero-vertices.

We distinguish this rectangular matrix from the full
Ferrers diagram $F(G)$. The latter is obtained from the
adjacency matrix in \eqref{eq:adj-Ferrers} by depicting
off-diagonal ones as black cells, off-diagonal zeros as white
cells, and marking the diagonal cells separately.
It is symmetric, and its strict lower-triangular part records
each edge exactly once, containing precisely $m$ black cells.

Conversely, a matrix of the block form
\eqref{eq:adj-Ferrers}, with a $(0,1)$ Ferrers matrix as its
off-diagonal block, defines a threshold graph.
This gives a matrix formulation of the correspondence between
creation sequences, nested neighborhoods, and Ferrers diagrams.

We now introduce two complementary tools used in this paper.
\subsection{The scalar comparison equation}

A common approach to spectral extremal problems is to compare the
coordinates of Perron vectors; see \cite{Bell1991,BrualdiHoffman1985,
BrualdiSolheid1986,ChangTam2011,CvetkovicRowlinson1988,
FanTamZhou2008,Friedland1985,Rowlinson1988,ZhaiLinZhao2022,ZhangJinWangYangZhang2025}.
 For threshold graphs, nested
neighborhoods facilitate qualitative comparisons, but the quantitative
estimates required to determine the effect of an edge rearrangement can
still be difficult to obtain.  
We use the standard reduction afforded by block symmetries and condensed 
matrices~\cite{FritscherTrevisan2016,HelmbergEtAl2025}, and work 
directly with the block Perron equations to construct an identity 
involving the Perron coordinates, the degrees, and the spectral radius.

Retain the notation introduced above. For \(\lambda\ge n\), put
\begin{equation}\label{eq:scalar-LR}
 L_i(\lambda)=\lambda-t_i,\qquad
 R_i(\lambda)=\lambda-b+2-s_i,\qquad
 R_0(\lambda)=\lambda-b+2,
\end{equation}
and define
\begin{equation}\label{eq:scalar-functions}
 S_G(\lambda)=\sum_{i=1}^k
       \frac{a_i}{L_i(\lambda)R_{i-1}(\lambda)R_i(\lambda)},\qquad
 T_{a,b}(\lambda)=\frac{1}{\lambda(\lambda-n+2)}
       -\frac{1}{2(\lambda-b+1)(\lambda-b+2)}.
\end{equation}
Since \(L_i(\lambda)\ge \lambda-b\ge a>0\) and
\(R_i(\lambda)\ge \lambda-n+2\ge2\), all denominators are positive when 
\(\lambda\in [n,\infty)\). In particular, \(S_G\) and \(T_{a,b}\) are continuous
on this interval.

\begin{proposition}\label{prop:kappa}
Let \(G\) be a connected threshold graph of order \(n\ge2\), with
vertex-sequence $(a_1,\ldots,a_k\,;\,b_1,\ldots,b_k),$
and let \(\lambda=\rho(G)\). Let  $\mathbf  x$ be a positive Perron vector of
\(Q(G)\), and denote its constant coordinates on the blocks \(A_i\)
and \(B_i\) by \(x_i\) and \(y_i\), respectively. Then
\[
\kappa:=L_i(\lambda)R_i(\lambda)(x_i+y_i)
\]
is independent of \(i\), for \(1\le i\le k\).
\end{proposition}
\begin{proof}
Every connected threshold graph has a dominating vertex and hence
contains a spanning star. Thus \(\lambda\ge n\), so all the quantities
in~\eqref{eq:scalar-LR} are positive. In the following,
we write \(L_i=L_i(\lambda)\) and \(R_i=R_i(\lambda)\).

The block eigenvalue equations are
\begin{equation}\label{eq:block-perron}
 L_i x_i=\sum_{j=i}^k b_jy_j,\qquad
 R_i y_i=\sum_{j=1}^k b_jy_j+\sum_{j=1}^i a_jx_j.
\end{equation}
For \(1\le i\le k-1\), consecutive equations in
\eqref{eq:block-perron} give
$$
 L_{i+1}(x_{i+1}+y_i)=L_i(x_i+y_i),~~
 R_{i+1}(x_{i+1}+y_{i+1})=R_i(x_{i+1}+y_i).
$$
Therefore, we have
\[
 L_{i+1}R_{i+1}(x_{i+1}+y_{i+1})
 =L_iR_i(x_i+y_i).
\]
Hence \(L_iR_i(x_i+y_i)\) is independent of \(i\), as proved.
\end{proof}

Using Proposition~\ref{prop:kappa}, we derive a one-point comparison criterion, which
 yields an exact scalar equation for the signless
Laplacian spectral radius.
 In particular, spectral improvement can be verified by evaluating a scalar 
 rational function, without estimating the individual Perron coordinates.
\begin{lemma}\label{lem:scalar-comparison}
Let \(G\) be a connected threshold graph of order \(n\ge2\), with
vertex-sequence as above, and let \(\lambda=\rho(G)\). Then
$
S_G(\lambda)=T_{a,b}(\lambda),
$
where \(a=\sum_{i=1}^k a_i\) and \(b=\sum_{i=1}^k b_i\).
Moreover, if \(H\) is a connected threshold graph with the same total
zero-block size \(a\) and total one-block size \(b\) as \(G\), and
\[
S_H(z)>T_{a,b}(z)
\]
for some \(z\ge n\), then \(\rho(H)>z\). The graphs \(G\) and \(H\)
need not have the same character.
\end{lemma}
\begin{proof}
The first assertion follows by telescoping the block Perron equations.
Let \(x_i,y_i\) and \(\kappa\) be as in
Proposition~\ref{prop:kappa}, and set
$
y_0=\frac{\sum_{i=1}^k b_i y_i}{R_0}.
$
Subtracting consecutive equations in the second part of
\eqref{eq:block-perron}, including the initial identity
\(R_0y_0=\sum_{i=1}^k b_i y_i\), gives
$
R_i y_i-R_{i-1}y_{i-1}=a_i x_i.
$
Since \(R_{i-1}-R_i=a_i\), Proposition~\ref{prop:kappa} implies
\[
y_i-y_{i-1}
 =\frac{a_i\kappa}{L_iR_{i-1}R_i},
 \qquad 1\le i\le k.
\]
At the first block, the equations in \eqref{eq:block-perron} give
\[
\kappa=(R_1+L_1+a_1)\sum_{i=1}^k b_i y_i
      =2(\lambda-b+1)\sum_{i=1}^k b_i y_i.
\]
At the last block, \(L_kx_k=b_ky_k\) and \(L_k+b_k=\lambda\), so
\(\kappa=\lambda R_k y_k\). Summing the preceding differences now yields
\[
S_G(\lambda)=\frac{y_k-y_0}{\kappa}
 =\frac{1}{\lambda R_k}
  -\frac{1}{2(\lambda-b+1)R_0}
 =T_{a,b}(\lambda),
\]
because \(R_k=\lambda-n+2\) and \(R_0=\lambda-b+2\).

For the comparison assertion, suppose first that $t\ge n$
satisfies $S_H(t)=T_{a,b}(t)$.
Use the same block notation for $H$, with $k$ now denoting
its character, and evaluate all $L_i,R_i$ at $t$.
Set $y_0=1/[2(t-b+1)R_0]$ and define
\[
y_i=y_0+\sum_{j=1}^i\frac{a_j}{L_jR_{j-1}R_j},~~
x_i=\frac{1}{L_iR_{i-1}}-y_{i-1}.
\]
All $y_i$ are positive. The definitions give
$x_i=1/(L_iR_i)-y_i$, while the scalar equation gives
$y_k=1/(tR_k)$. Consequently,
\[
x_k=\frac{b_k}{tL_kR_k}>0,\qquad L_kx_k=b_ky_k.
\]
Moreover, $L_ix_i=L_{i+1}x_{i+1}+b_iy_i$ for $i<k$.
Backward induction therefore yields
$L_ix_i=\sum_{j=i}^k b_jy_j>0$ for every $i$. In particular,
\[
\sum_{j=1}^k b_jy_j=L_1x_1=\frac1{R_0}-L_1y_0=R_0y_0.
\]
The definitions also give
$R_iy_i-R_{i-1}y_{i-1}=a_ix_i$, and hence
\[
R_iy_i=\sum_{j=1}^k b_jy_j+\sum_{j=1}^i a_jx_j.
\]
Thus these positive coordinates satisfy the block eigenvalue
equations of $Q(H)$. By the Perron--Frobenius theorem,
$t=\rho(H)$. Together with the first part applied to $H$,
this shows that $S_H(t)=T_{a,b}(t)$ has the unique solution
$t=\rho(H)$ on $[n,\infty)$.

Now $S_H(t)-T_{a,b}(t)$ is continuous on $[n,\infty)$, and
\[
\lim_{t\to\infty}t^2\bigl(S_H(t)-T_{a,b}(t)\bigr)=-\frac12.
\]
Hence this difference is negative for all sufficiently large $t$.
Since $S_H(z)>T_{a,b}(z)$, the intermediate value theorem
provides a zero strictly larger than $z$.
By the preceding characterization, this zero is $\rho(H)$.
Therefore $\rho(H)>z$.
\end{proof}

\subsection{$Q$-Coronal for threshold graphs}
\label{sec:coronal}

In~\cite{McLemanMcNicholas2011}, McLeman and McNicholas introduced the coronal of a
graph \(G\) of order \(n\), defined as the sum of the entries of
\((\lambda I_n-A)^{-1}\), where \(A\) is the adjacency matrix of \(G\).
This concept was generalized to a graph matrix \(M\) of \(G\) by Cui and
Tian~\cite{CuiTian2012}. In this paper, we take \(M=Q(H)\), the signless Laplacian matrix
of a graph \(H\), as follows.

\begin{definition}
Let $G$ be a graph of order $n$, and put
$\phi_G(\lambda)=\det(\lambda I_n-Q(G))$.
Its $Q$-coronal is the rational function
\[
\Gamma_G(\lambda)=\mathbf1_n^{\mathsf T}(\lambda I_n-Q(G))^{-1}\mathbf1_n.
\]
The inverse is taken over $\mathbb C(\lambda)$, where it exists
because $\phi_G(\lambda)$ is a monic polynomial.
\end{definition}

The \(Q\)-coronal records the resolvent of \(Q(G)\) in the all-ones
direction. Together with \(\phi_G(\lambda)\), it
provides the data needed to track spectral information under threshold
graphs. We first record how both quantities change when a vertex is
appended to the creation word.

If a symbol \(0\) is appended, the new graph is \(H\sqcup K_1\), and
\begin{equation}\label{eq:threshold-isolated-extension}
\phi_{H\sqcup K_1}(\lambda)=\lambda\phi_H(\lambda),
\qquad
\Gamma_{H\sqcup K_1}(\lambda)
 =\Gamma_H(\lambda)+\frac1\lambda.
\end{equation}
If a symbol \(1\) is appended, the new graph is \(K_1\vee H\). With
the new dominating vertex listed first, its signless Laplacian matrix is
\[
Q(K_1\vee H)=
\begin{pmatrix}
n & \mathbf{1}_n^{\mathsf T}\\
\mathbf{1}_n & Q(H)+I_n
\end{pmatrix}.
\]
Taking the Schur complement of \((\lambda-1)I_n-Q(H)\) gives
\begin{equation}\label{eq:threshold-dominating-characteristic}
\phi_{K_1\vee H}(\lambda)
 =\phi_H(\lambda-1)
  \bigl(\lambda-n-\Gamma_H(\lambda-1)\bigr).
\end{equation}
Moreover, block inversion gives
\begin{equation}\label{eq:threshold-dominating-coronal}
\Gamma_{K_1\vee H}(\lambda)
 =\Gamma_H(\lambda-1)
  +\frac{\bigl(1+\Gamma_H(\lambda-1)\bigr)^2}
         {\lambda-n-\Gamma_H(\lambda-1)}.
\end{equation}
These identities are understood over \(\mathbb{C}(\lambda)\).

We now apply these extension formulas to a pair of threshold graphs.
For a graph \(F\), let \(\phi_F\) and \(\Gamma_F\) denote its
\(Q\)-characteristic polynomial and \(Q\)-coronal, respectively.

\begin{theorem}\label{thm:threshold-cospectral-transform}
Let \(k\ge2\), and let \(a_2,\ldots,a_k\) and \(b_2,\ldots,b_k\)
be positive integers. The threshold graphs
$
G=(3,a_2,\ldots,a_k\,;\,1,b_2,\ldots,b_k)$ and
$G'=(1,a_2+1,a_3,\ldots,a_k\,;\,2,b_2,\ldots,b_k)
$
are \(Q\)-cospectral. In particular, \(\rho(G)=\rho(G')\).
\end{theorem}

\begin{proof}
The corresponding binary creation words are
\[
0001\,0^{a_2}1^{b_2}\cdots 0^{a_k}1^{b_k}
\quad\text{and}\quad
0110\,0^{a_2}1^{b_2}\cdots 0^{a_k}1^{b_k}.
\]
Thus the graphs have cores \(K_{1,3}\) and \(K_3\sqcup K_1\),
respectively, followed by the same sequence of threshold extensions.
A direct computation shows that the two cores have the same \(Q\)-characteristic polynomial and \(Q\)-coronal:
\[
\begin{aligned}
\phi_{K_{1,3}}(\lambda)
&=\phi_{K_3\sqcup K_1}(\lambda)
=\lambda(\lambda-1)^2(\lambda-4), ~~\text{and}~~
\Gamma_{K_{1,3}}(\lambda)
&=\Gamma_{K_3\sqcup K_1}(\lambda)
=\frac1\lambda+\frac3{\lambda-4}.
\end{aligned}
\]

Starting from these cores, append the common suffix
\(0^{a_2}1^{b_2}\cdots 0^{a_k}1^{b_k}\) one symbol at a time.
Equations~\eqref{eq:threshold-isolated-extension}--%
\eqref{eq:threshold-dominating-coronal} show inductively that after
each common extension the two graphs still have equal
\(Q\)-characteristic polynomials and equal \(Q\)-coronals. In
particular, their final \(Q\)-characteristic polynomials are equal.
Therefore their \(Q\)-spectra, including multiplicities, agree, so
\(G\) and \(G'\) are \(Q\)-cospectral.
\end{proof}

For the reduction below, a threshold graph with vertex-sequence
$(a_1,\ldots,a_k;b_1,\ldots,b_k)$ is called \emph{nonregular} if
$k\ge2$ and $(a_1,b_1)=(1,2)$; otherwise it is called \emph{regular}.
As a consequence of Theorem \ref{thm:threshold-cospectral-transform}, we have

\begin{corollary}\label{cor:pi11-pi13}
Let $n\ge5$. Every nonregular threshold graph
$
  G=(1,a_2,\ldots,a_k;\,2,b_2,b_3,\ldots,b_k)
$
is $Q$-cospectral with the regular threshold graph $\widehat G$
obtained by replacing the initial data
\[
  (a_1,b_1,a_2)=(1,2,a_2)
  \quad\text{with}\quad
  (3,1,a_2-1).
\]
If $a_2=1$, the resulting zero block is deleted and the adjacent
one-blocks are merged. Consequently, every nonregular
$Q$-maximal threshold graph has a regular $Q$-maximal
representative with the same $Q$-spectrum.

Conversely, let
$
  \widehat G=(a_1,\ldots,a_\ell;\,b_1,\ldots,b_\ell)
$
be a regular $Q$-maximal threshold graph. If $a_1=3$, an additional
nonregular $Q$-maximal graph is obtained by the inverse
transformation:
\[
  G=
  \begin{cases}
    (1,a_2+1,a_3,\ldots,a_\ell;\,
       2,b_2,b_3,\ldots,b_\ell),
       & b_1=1,\\[2mm]
    (1,1,a_2,\ldots,a_\ell;\,
       2,b_1-1,b_2,\ldots,b_\ell),
       & b_1\ge2.
  \end{cases}
\]
Here, empty tails are omitted. These inverse transformations
recover all nonregular $Q$-maximal threshold graphs.
\end{corollary}

Thus, once the regular $Q$-maximal graphs have been determined,
all nonregular $Q$-maximal graphs are recovered from those regular
representatives whose initial zero-block has size three.
In particular, when $j'=3$, the graph $S^\triangle_{n,m}$ is
$Q$-cospectral with $S_{n,m}$ and is therefore $Q$-maximal whenever
$S_{n,m}$ is $Q$-maximal.

\section{Structural properties and spectral estimates}
\label{sec:qmax-structure}

To streamline the exposition and preserve the flow of the main
argument, we state the structural properties and spectral estimates
in this section and collect their proofs in
Section~\ref{sec:proofs-structure}.
%{\color{red}In view of Corollary~\ref{cor:pi11-pi13}, our structural analysis
%focuses on regular $Q$-maximal threshold graphs.}

To analyze local rearrangements, we first rewrite $ S_G(z)$ in
\eqref{eq:scalar-functions} as a sum over individual zero positions.
For $z\ge n$, define
\[
\begin{aligned}
  r_\ell(z)&=z-b+2-\ell,~~
  \omega_\ell(z)=
    \frac{1}{r_{\ell-1}(z)r_\ell(z)},
    &&1\le\ell\le a.
\end{aligned}
\]

Since $r_{s_i}(z)=R_i(z)$ and
$r_{\ell-1}(z)-r_\ell(z)=1$, telescoping gives
\[
  \sum_{\ell=s_{i-1}+1}^{s_i}\omega_\ell(z)
  =\frac{1}{R_i(z)}-\frac{1}{R_{i-1}(z)}
  =\frac{a_i}{R_{i-1}(z)R_i(z)}.
\]
Consequently,
\begin{equation}\label{eq:positional-S}
  S_G(z)
  =\sum_{i=1}^k\frac{1}{L_i(z)}
     \sum_{\ell=s_{i-1}+1}^{s_i}\omega_\ell(z).
\end{equation}
For $z\ge n$, the quantities $r_\ell(z)$ are positive and
strictly decreasing with $\ell$. Thus the positional weights
$\omega_\ell(z)$ are positive and strictly increasing with $\ell$.

The following lemma gives three local restrictions on
consecutive block sizes.

\begin{lemma}\label{lem:qmax-blocks}
Let $G\in\mathcal{C}_{n,m}$ be a regular $Q$-maximal threshold graph with
vertex-sequence
$(a_1,\ldots,a_k;\,b_1,\ldots,b_k).$
Then the following statements hold:
\begin{enumerate}[label=\textup{(\roman*)}]
\item
For every $1\le i\le k-1$, either $b_i=1$ or $
  b_i>a_i+a_{i+1}.$
\item
If $1\le i\le k-2$ and $b_i=1$, then $
  b_{i+1}+1>a_i+a_{i+1}+a_{i+2}.$
\item
If $1\le i\le k-2$, $b_i>1$, and $b_{i+1}=1$, then
$
  b_i>a_{i+1}+a_{i+2}.
$
\end{enumerate}
\end{lemma}

\begin{remark}
Consider a Ferrers
 diagram of a (regular) threshold graph, and order its vertices by
$ (B_k,\ldots,B_1,A_1,\ldots,A_k).$
We restrict attention to its lower-triangular part. Number the rows and columns from $1$ to $n$,
and call $r+c$ the \emph{coordinate sum} of an element $(r,c)$,
where $r>c$.

The local exchanges used in the proof of
Lemma~\ref{lem:qmax-blocks} exhibit a directional asymmetry.
For the downward moves covered by that proof, a non-decrease
in coordinate sum yields a strict increase in the signless
Laplacian spectral radius. This implication cannot simply be
extended to upward moves: even a strict increase in coordinate
sum need not increase the spectral radius.

For example, let $G=(2,3;\,3,2)$ and
$H=(1,2,1;\,4,1,1)$, both belonging to $\mathcal C_{10,26}$.
Let $
  \mathcal F(G)
  =\{(r,c):1\le c<r\le n,\ F(G)_{rc}=1\}
$
denote the support of the lower-triangular part of the Ferrers diagram, 
and define $\mathcal F(H)$ analogously.
Then
$
  \mathcal F(H)
  =\bigl(\mathcal F(G)\setminus\{(10,2)\}\bigr)
     \cup\{(7,6)\}.
$
Thus a single entry moves upward, while its coordinate sum
increases from $12$ to $13$. Nevertheless,
$ \rho(G)>13.14
  >\rho(H).$
\end{remark}

\begin{lemma}
\label{lem:qmax-singletons}
Let $G\in\mathcal C_{n,m}$ be a $Q$-maximal threshold graph
with vertex-sequence $(a_1,\ldots,a_k;$ $\,b_1,\ldots,b_k)$.
Then the following statements hold:
\begin{enumerate}[label=\textup{(\roman*)}]
\item
There is at most one index $i\in\{1,\ldots,k-1\}$ such that
$b_i=1$.

\item
If $a_i=b_j=1$ for some $ 2\le i\le k-1$ and 
  $1\le j\le k-1,$
then $j\in\{i-1,i\}$.
\end{enumerate}
\end{lemma}

The next lemma gives a strict spectral improvement under
an admissible symmetric perturbation of the one-block sizes.

\begin{lemma}\label{lem:strict-midpoint}
Let $G\in\mathcal C_{n,m}$ be a threshold graph with vertex-sequence
$(a_1,\ldots,a_k;\,b_1,\ldots,b_k)$.
Suppose that a nonzero integer vector $(\delta_i)_{i=1}^k$
satisfies $\sum_{i=1}^k\delta_i=0$ and
$\sum_{i=1}^k s_i\delta_i=0$.
Assume also that $|\delta_i|\le b_i$ for $1\le i<k$
and $|\delta_k|<b_k$.
Let $G_\pm=(a_1,\ldots,a_k;\,b_1\pm\delta_1,\ldots,b_k\pm\delta_k)$.
Any internal one-block of size zero is suppressed,
and the adjacent zero-blocks are merged.
Then
\[
\rho(G)<\max\{\rho(G_+),\rho(G_-)\}.
\]
Consequently, no $Q$-maximal threshold graph admits such a
nonzero perturbation.
\end{lemma}

We next consider a local operation that splits a zero-block
into two nonempty parts and redistributes the neighbouring
one-blocks while preserving both the order and the size.

Fix $2\le i\le k-1$ and $\tau\in\{2,3\}$.
An integer $r$ is called admissible if
\begin{equation}\label{eq:split-range}
  \max\left\{
    1,a_i-\left\lfloor\frac{b_{i-1}}{\tau}\right\rfloor
  \right\}
  \le r\le
  \min\left\{
    a_i-1,\left\lfloor\frac{b_i}{\tau}\right\rfloor
  \right\}.
\end{equation}
For such an $r$, replace the local block data by
\begin{equation*}
\begin{aligned}
  (a_{i-1},a_i;\,b_{i-1},b_i)
  \longmapsto
  \bigl(a_{i-1},r,a_i-r;\,
    b_{i-1}-\tau(a_i-r),\tau a_i,b_i-\tau r\bigr),
\end{aligned}
\end{equation*}
leaving all other block sizes unchanged.
Denote the resulting graph by $G'$.

Condition~\eqref{eq:split-range} ensures that the two new
zero-blocks have positive sizes and that all resulting
one-block sizes are nonnegative.
Any block of size zero is deleted, and adjacent blocks of
the same type are merged.
The total zero-block size and the total one-block size
are unchanged, and a direct edge count gives $
  |V(G')|=|V(G)|$ and $
  |E(G')|=|E(G)|.$
Moreover, the final one-block $B_k$ is unchanged, so $G'$
remains connected. Thus $G'\in\mathcal{C}_{n,m}$ whenever $G\in\mathcal{C}_{n,m}$.

\begin{lemma}
\label{lem:splitting-inequality}
Let $G\in \mathcal{C}_{n,m}$ be a $Q$-maximal threshold graph with
vertex-sequence $
  (a_1,\ldots,a_k;$ $\,b_1,\ldots,b_k).$
For every $2\le i\le k-1$, every $\tau\in\{2,3\}$,
and every integer $r$ satisfying \eqref{eq:split-range},
there has
\begin{equation}\label{eq:splitting-bound}
  \tau\bigl[R_i(\rho(G))+a_i-r\bigr]
  -L_i(\rho(G))\le0.
\end{equation}
In particular, if $a_i\ge2$ and $
  b_{i-1}\ge
    2\left\lceil\frac{a_i}{2}\right\rceil,
  b_i\ge
    2\left\lfloor\frac{a_i}{2}\right\rfloor,
$
then the choice
$\tau=2$ and $r=\lfloor a_i/2\rfloor$ is admissible, and
\begin{equation}\label{eq:splitting-special-bound}
\begin{aligned}
  \rho(G)
  \le
  \sum_{j=1}^k b_j+\sum_{j=1}^{i-1}b_j
  +2\sum_{j=1}^{i-1}a_j
  +2\left\lfloor\frac{a_i}{2}\right\rfloor-4.
\end{aligned}
\end{equation}
\end{lemma}

For the second and third zero-blocks, the splitting
inequality, together with the clique-join and prefix
estimates below, gives the following restrictions.

\begin{lemma}
\label{lem:splitting-consequences}
Let $G\in\mathcal{C}_{n,m}$ be a regular $Q$-maximal threshold graph with
vertex-sequence $
  (a_1,\ldots,a_k;\,b_1,\ldots,b_k).$
Suppose that $i\in\{2,3\}$ satisfies $i\le k-1$,
$a_i\ge2$, and $b_{i-1},b_i>1$.
Then
\[
  b_{i-1}+b_i\le3a_i+1,
  \qquad
  a_i\ge a_{i-1}+a_{i+1}+1.
\]
\end{lemma}

We next record two auxiliary spectral estimates.
The first describes the spectral growth caused by joining a clique,
while the second provides lower bounds in terms of the numbers of
edges between zero-blocks and one-blocks.

\begin{lemma}
\label{lem:clique-join-estimate}
Let $G$ be a graph of order $n\ge1$.
Then, for every integer $t\ge1$,
\[
  \rho(G\join K_t)\ge \rho(G)+2t.
\]
\end{lemma}

\begin{lemma}\label{lem:rayleigh-lower-bounds}
Let $G\in\mathcal C_{n,m}$ be a threshold graph with
vertex-sequence $(a_1,\ldots,a_k;\,b_1,\ldots,b_k)$,
and put $a=\sum_{j=1}^k a_j$, $b=\sum_{j=1}^k b_j$,
and $n=a+b$.
Then
\[
\rho(G)>\frac32n-2
\]
holds if at least one of the following conditions is satisfied:
\begin{enumerate}[label=\textup{(\arabic*)}]
\item
The total number of edges between zero-blocks and one-blocks
satisfies
$\sum_{j=1}^k b_js_j\ge a^2.$

\item
There exists an index $1\le\ell<k$ such that
$\sum_{j=1}^{\ell}b_js_j\ge s_\ell^2$ and $\sum_{j=\ell+1}^k b_j(s_j-s_\ell)\ge(a-s_\ell)^2.$
\end{enumerate}
\end{lemma}

The following lemma applies the preceding Rayleigh bounds to
the initial one-pair and two-pair subgraphs of a regular
$Q$-maximal threshold graph.

\begin{lemma}
\label{lem:prefix-estimates}
Let $G$ be a regular $Q$-maximal threshold graph with
vertex-sequence $
  G=(a_1,\ldots,a_k;\,b_1,\ldots,b_k), k\ge3.$
For $r\in\{1,2\}$, let $H_r$ be the subgraph of $G$ induced
by the first $r$ zero--one block pairs. 
Then the following statements hold:
\begin{enumerate}[label=\textup{(\arabic*)}]
\item
If $b_1>1$, then $
  \rho(H_1)>\frac32(a_1+b_1)-2.$
\item
If $b_2>1$, then $
  \rho(H_2)>
  \frac32(a_1+b_1+a_2+b_2)-2.$
\end{enumerate}
\end{lemma}

We also require two parameter-dependent lower bounds for
the short threshold graphs occurring in the subsequent
compression arguments.

\begin{lemma}
\label{lem:revised-prefix-bounds}
The following statements hold:
\begin{enumerate}[label=\textup{(\arabic*)}]
\item
Let $F$ be a threshold graph with vertex-sequence $
  F=(a_1,a_2;\,b_1,b_2).
$
If $
  b_1\ge a_1+a_2-1,
  b_2\ge a_2+2,
  b_1+b_2\le3a_2+2,$
then
\[
  \rho(F)>2b_1+b_2+2a_1+a_2-4.
\]
\item
Let $F$ be a threshold graph with vertex-sequence $
  F=(a_1,1,a_3;\,b_1,b_2,b_3).$
If $
    a_1+2\le b_1\le a_3-1,
  a_3+2\le b_2\le2a_3,
    b_3\ge a_3+2,$
then
\[
  \rho(F)>2b_1+2b_2+b_3+2a_1+a_3-3.
\]
\end{enumerate}
\end{lemma}

\section{Proof of theorems \ref{thm:main1} and \ref{thm:main2}}

\subsection{Proof of theorem \ref{thm:main1}}

\begin{lemma}
\label{lem:singleton-concentration}
Let $G\in\mathcal C_{n,m}$ be a regular $Q$-maximal
threshold graph with vertex-sequence
$(a_1,\ldots,a_k;\,b_1,\ldots,b_k), k\ge4.$
Then $b_2>1$.
\end{lemma}

\begin{proof}
Suppose that $b_2=1$. By Lemmas~\ref{lem:qmax-singletons}\textup{(i)} and~\ref{lem:qmax-blocks}\textup{(i)--(iii)},
\[
b_1\ge\max\{a_1+a_2+1,a_2+a_3+1\},\qquad b_3\ge a_2+a_3+a_4.
\]
In particular, $L_4=L_1+b_1+b_3+1$ and $R_3=R_1-a_2-a_3$.
The preceding bounds imply
\[
b_1+b_3+1\ge2(a_2+a_3)+3,\qquad a_1+a_2+a_3\le b_1+b_3-3.
\]
The prefix formed by the first three block pairs satisfies $a_1b_1\ge a_1^2$ and $a_2(b_3+1)+a_3b_3>(a_2+a_3)^2$.
Its join with $K_{t_4}$ is a subgraph of $G$, so Lemmas~\ref{lem:rayleigh-lower-bounds}\textup{(2)} and~\ref{lem:clique-join-estimate} give
\[
\lambda>\frac32(a_1+a_2+a_3+b_1+b_3+1)+2t_4-2.
\]
Also, $G$ properly contains $K_{b+1}$, whence $\lambda>2b$ and $2L_1>L_4$.
Using $L_1=R_1+a_1-2$ and the preceding bounds, we obtain
\[
3R_3>L_4+4,\qquad \frac{L_4}{R_1}<\frac{3b_1+3b_3-4}{a_2+a_3+b_1+b_3}<\frac{b_1+b_3+1}{a_2+a_3}.
\]

Choose integers $\ell,\varepsilon$ such that
\[
(a_2+a_3)b_3+a_2=\ell(b_1+b_3+1)+\varepsilon,\qquad 0\le\varepsilon<b_1+b_3+1.
\]
Then $0\le\ell\le a_2+a_3-1$.
Replace the first four block pairs of $G$ by
\[
(a_1+\ell,1,a_2+a_3+a_4-\ell-1;\,b_1+b_3+1-\varepsilon,\varepsilon,b_4),
\]
leaving the remaining blocks unchanged, and call the resulting graph $G'$.
If $\varepsilon=0$, merge the adjacent zero-blocks.
All remaining block sizes are positive, the totals $a,b$ are preserved, and
$e(G')-e(G)=\ell(b_1+b_3+1)+\varepsilon-[(a_2+a_3)b_3+a_2]=0.$
Hence $G'\in\mathcal C_{n,m}$.
Define
\[
D=\frac{\ell(b_1+b_3+1)}{L_1(R_1-\ell)}+\frac{\varepsilon R_1}{(L_4-\varepsilon)(R_1-\ell)(R_1-\ell-1)}-\frac{(a_2+a_3)[a_2(b_3+1)+a_3b_3]}{(a_2L_2+a_3L_3)R_3}.
\]
All denominators are positive, since $R_3\ge a_4+2$ and $L_4-\varepsilon>L_1>0$.
We have
\[
\frac{a_2(b_3+1)}{L_2R_1R_2}+\frac{a_3b_3}{L_3R_2R_3}<\frac{(a_2+a_3)[a_2(b_3+1)+a_3b_3]}{(a_2L_2+a_3L_3)R_1R_3}.
\]
Indeed, clearing positive denominators reduces this to $(a_2+a_3)L_3>R_1$, which follows from $L_3-R_1=a_1+b_1-1>0$.
Substitution into the scalar formula and cancellation of the unchanged terms yield
$L_4R_1\bigl(S_{G'}(\lambda)-S_G(\lambda)\bigr)>D.$
It remains to prove $D>0$.

\smallskip
\noindent\textbf{Case 1. $\varepsilon=0$.}
Then $\ell>0$. By the division identity, $D>0$ is equivalent to
$(b_1+b_3+1)R_3-(a_2+a_3)L_1>0.$
This follows from the upper bound for $L_4/R_1$, since
$
(b_1+b_3+1)R_3-(a_2+a_3)L_1=(b_1+b_3+1)R_1-(a_2+a_3)L_4>0.$

\smallskip
\noindent\textbf{Case 2. $\varepsilon>0$ and $\ell=0$.}
By the division identity, $D>0$ is equivalent to
$[(a_2+a_3)b_3+a_2]R_1-(a_2+a_3)L_4>0.$
Since $b_3\ge3$ and $3R_3>L_4+4$, we have $b_3R_1>3R_3>L_4$. Hence
$[(a_2+a_3)b_3+a_2]R_1-(a_2+a_3)L_4=(a_2+a_3)(b_3R_1-L_4)+a_2R_1>0.$

\smallskip
\noindent\textbf{Case 3. $\varepsilon>0$ and $\ell=a_2+a_3-1$.}
By the division identity, $D>0$ is equivalent to
$[(a_2+a_3)b_1+a_3]R_3-(a_2+a_3)L_1>0.$
Since $b_1\ge3$, we have $b_1R_3\ge3R_3>L_4>L_1$. Hence
$[(a_2+a_3)b_1+a_3]R_3-(a_2+a_3)L_1=(a_2+a_3)(b_1R_3-L_1)+a_3R_3>0.$

\smallskip
\noindent\textbf{Case 4. $\varepsilon>0$ and $1\le\ell\le a_2+a_3-2$.}
In particular, $a_2+a_3\ge3$.
Put
\[
p=\frac{(a_2+a_3)(b_1+b_3+1)}{\min\{a_2b_1+a_3(b_1+1),\,a_2(b_3+1)+a_3b_3\}}.
\]
The division identity and the bounds on the block sizes give
$2\le p<a_2+a_3$ and $b_1+b_3+1\ge(a_2+a_3+1)p.$
For $s\ge2$ and $z\ge2s+2$, define
\[
f(s,z)=\frac{(s+1)(z^2-(2s+1)z+7s)}{z(z+s-1)}.
\]
The upper bound for $L_4/R_1$ gives
\[
1-\frac{(a_2+a_3)L_4}{(b_1+b_3+1)R_1}>\frac{f(a_2+a_3,b_1+b_3+1)}{a_2+a_3+1}.
\]
Moreover, AM--GM and $2L_1>L_4$ imply
\[
\frac{\varepsilon(b_1+b_3+1-\varepsilon)}{L_4-\varepsilon}\le(\sqrt{L_4}-\sqrt{L_1})^2<\frac{3(L_4-L_1)^2}{7L_1+5L_4}.
\]
Since $\ell\le a_2+a_3-2$ and $0<\varepsilon<b_1+b_3+1$, we have
$R_1-\ell-\frac{\varepsilon}{b_1+b_3+1}>R_3+1>0.$
Multiply $D$ by the positive factor
$\frac{L_1R_3}{R_1}\left(R_1-\ell-\frac{\varepsilon}{b_1+b_3+1}\right).$
Using the division identity and
\[
\frac{R_3}{R_1-\ell}\left(1-\frac{L_1}{(b_1+b_3+1)(R_1-\ell-1)}\right)<1,
\]
these estimates show that the following inequality is sufficient for $D>0$:
\[
\frac{(7L_1+5L_4)(p-1)(a_2+a_3)^2}{3p^2(a_2L_2+a_3L_3)(a_2+a_3+1)}f(a_2+a_3,b_1+b_3+1)>1.
\]
To verify it, the definition of $p$ and $2L_1>L_4$ give
\[
\frac{(7L_1+5L_4)(a_2+a_3)}{a_2L_2+a_3L_3}\ge\frac{7L_1+5L_4}{L_4-(L_4-L_1)/p}>\frac{17p}{2p-1}.
\]
Also, $f>0$, and direct differentiation shows that $f$ increases in $z$ and that $s\mapsto f(s,p(s+1))$ increases for $s\ge p$.
Since $a_2+a_3>p$ and $b_1+b_3+1\ge p(a_2+a_3+1)$,
the preceding monotonicity properties give
\[
f(a_2+a_3,b_1+b_3+1)\ge f(p,p(p+1))=\frac{p^3-2p+6}{p^2+2p-1}.
\]
Moreover, $a_2+a_3\ge3$ gives
$\frac{a_2+a_3}{a_2+a_3+1}\ge\frac34.$
The bounds on $p$ and $a_2+a_3\ge3$ therefore yield
\[
\frac{a_2+a_3}{a_2+a_3+1}f(a_2+a_3,b_1+b_3+1)\ge\frac{3(p^3-2p+6)}{4(p^2+2p-1)}.
\]
The desired inequality follows by multiplying these two bounds, since
\[
\frac{17(p-1)(p^3-2p+6)}{4p(2p-1)(p^2+2p-1)}=1+\frac{2+(p-2)^2(9p^2+7p-26)}{4p(2p-1)(p^2+2p-1)}>1.
\]
Thus $D>0$ in every case. Lemma~\ref{lem:scalar-comparison} now gives $\rho(G')>\rho(G)$, contradicting the $Q$-maximality of $G$.
Therefore $b_2>1$.
\end{proof}

\begin{proof}[Proof of Theorem~\ref{thm:main1}]
By Lemma~\ref{lemma1}, every $Q$-maximal graph in $\mathcal{C}_{n,m}$
is a threshold graph.
We first exclude regular $Q$-maximal threshold graphs of
character at least four, and then treat the nonregular
case.

Throughout the proof, put $\lambda=\rho(G)$ and retain
the notation $
  t_j=\sum_{h=j}^k b_h.$
We repeatedly use the following observation.
If $F$ is a subgraph of the prefix through $B_r$, then
$F\join K_{t_{r+1}}$ is a subgraph of $G$.
Consequently, Lemma~\ref{lem:clique-join-estimate} gives
\begin{equation}\label{eq:prefix-clique-comparison}
  \lambda
  \ge \rho(F\join K_{t_{r+1}})
  \ge \rho(F)+2t_{r+1}.
\end{equation}
These clique joins are used only to obtain lower bounds;
they are not competitors required to belong to $\mathcal{C}_{n,m}$.

For an internal zero-block $A_i$, write
\begin{equation*}
  U_i=
  \sum_{j=1}^k b_j+\sum_{j=1}^{i-1}b_j
  +2\sum_{j=1}^{i-1}a_j
  +2\left\lfloor\frac{a_i}{2}\right\rfloor-4.
\end{equation*}
Whenever the balanced split with $\tau=2$ is admissible,
Lemma~\ref{lem:splitting-inequality} implies
$ \lambda\le U_i.$
Assume to the contrary that $
  G=(a_1,\ldots,a_k;\,b_1,\ldots,b_k), k\ge4,
$
is a regular  $Q$-maximal. We  show that no such graph exists by considering the following three cases.

\noindent {\bf Case 1. $b_1=1$.}
Lemma~\ref{lem:qmax-singletons}\textup{(i)} gives
$b_2,b_3>1$.
Moreover, $a_3=1$ would contradict
Lemma~\ref{lem:qmax-singletons}\textup{(ii)}, applied to
$A_3$ and $B_1$.
Thus $a_3\ge2$.
The block inequalities and
Lemma~\ref{lem:splitting-consequences} yield $
  b_2\ge a_1+a_2+a_3,
  b_3\ge a_3+2$ and $
  b_2+b_3\le3a_3+1.$
Delete the edges between $A_1$ and $B_1$, and restrict
the resulting graph to the prefix through $B_3$.
The resulting threshold graph is $
  F=(a_1+a_2+1,a_3;\,b_2,b_3).$
Its parameters satisfy the hypotheses of
Lemma~\ref{lem:revised-prefix-bounds}\textup{(1)}.
Hence \eqref{eq:prefix-clique-comparison} gives $
  \lambda>
  2b_2+b_3+2a_1+2a_2+a_3-2+2t_4.$

On the other hand, the displayed lower bounds on $b_2$
and $b_3$ make the balanced split of $A_3$ admissible.
Therefore, $
  \lambda\le U_3 \le
  2b_2+b_3+2a_1+2a_2+a_3-2+t_4.$
Since $t_4\ge1$, these inequalities are incompatible.

\noindent {\bf Case 2. $b_1>1$ and $b_2>1$.}
Suppose first that $a_2\ge2$.
The block inequalities and
Lemma~\ref{lem:splitting-consequences} give $
  b_1\ge a_1+a_2+1,
  b_2\ge a_2+2,
  b_1+b_2\le3a_2+1.$
Applying
Lemma~\ref{lem:revised-prefix-bounds}\textup{(1)}
to the first two block pairs and then using
\eqref{eq:prefix-clique-comparison}, we obtain $
  \lambda>
  2b_1+b_2+2a_1+a_2-4+2t_3.$
The balanced split of $A_2$ is admissible, whereas $
  \lambda\le U_2
  \le2b_1+b_2+2a_1+a_2-4+t_3.$
This is a contradiction. Therefore $a_2=1$.

If $b_3=1$, the pair $A_2,B_3$ contradicts
Lemma~\ref{lem:qmax-singletons}\textup{(ii)}.
Thus $b_3>1$.
We next show that $b_1<a_3$.
Otherwise, consider the symmetric perturbations
\[
  (b_1,b_2,b_3)
  \longmapsto
  \bigl(b_1\pm a_3,\,
        b_2\mp(a_3+1),\,
        b_3\pm1\bigr),
\]
with all other block sizes unchanged.
Both directions are admissible because $
  b_1\ge a_3,
  b_2\ge a_3+2,
  b_3>1.$
The total perturbation is zero, and $
  a_3s_1-(a_3+1)s_2+s_3=0,$
since $a_2=1$.
Lemma~\ref{lem:strict-midpoint} therefore gives a
contradiction.
Consequently, $
  b_1\le a_3-1.$
In particular, $a_3\ge2$.

The block inequalities now give $
  a_1+2\le b_1\le a_3-1,
  b_2\ge a_3+2,
  b_3\ge a_3+2.$
Furthermore,
Lemma~\ref{lem:splitting-consequences}, applied at $i=3$,
gives $
  b_2+b_3\le3a_3+1,$
and hence $b_2\le2a_3-1$.
The hypotheses of
Lemma~\ref{lem:revised-prefix-bounds}\textup{(2)}
are therefore satisfied.
Using \eqref{eq:prefix-clique-comparison}, we obtain
\[
  \lambda>
  2b_1+2b_2+b_3+2a_1+a_3-3+2t_4.
\]
The balanced split of $A_3$ is admissible, so
\[
  \lambda\le U_3
  \le
  2b_1+2b_2+b_3+2a_1+a_3-2+t_4.
\]
Since $t_4\ge1$, the preceding strict lower bound
contradicts this upper bound.
Notice that equality between the two comparison bounds
when $t_4=1$ causes no difficulty, because the lower
bound for $\lambda$ is strict.

\noindent {\bf Case 3.  $b_2=1$.}
By Lemma~\ref{lem:singleton-concentration}, there is no \(Q\)-maximal graph with character at least four 
and \(b_2=1\). Therefore, such a case can not occur.

The three cases exhaust all possibilities.
Thus every regular $Q$-maximal threshold graph has
character at most three.

It remains to exclude a $Q$-maximal threshold graph
\[
  G=(1,a_2,a_3,\ldots,a_k;\,
     2,b_2,b_3,\ldots,b_k),
  \qquad k\ge4.
\]
By Corollary~\ref{cor:pi11-pi13}, $G$ has a regular
$Q$-cospectral threshold graph representative $\widehat G$ , which is
also $Q$-maximal.

If $a_2\ge2$, then $
  \widehat G
  =(3,a_2-1,a_3,\ldots,a_k;\,
    1,b_2,b_3,\ldots,b_k).$
Its character is $k\ge4$, contradicting the regular case.

If $a_2=1$, normalization gives $
  \widehat G
  =(3,a_3,\ldots,a_k;\,
    b_2+1,b_3,\ldots,b_k).$
Its character is $k-1$.
The regular case therefore forces $k=4$.
However, $\widehat G$ then has three block pairs,
its first zero-block has size three, and its first
one-block has size $b_2+1\ge2$.
This contradicts
Lemma~\ref{lem:regular-three-pair-boundary}.

Thus no nonregular $Q$-maximal threshold graph has
character at least four.
Together with the regular case, this completes the
proof of Theorem~\ref{thm:main1}.
\end{proof}

\subsection{Proof of theorem \ref{thm:main2}}

\begin{lemma}\label{lem:residual-bound}
Let $G\in \mathcal{C}_{n,m}$ be a regular  $Q$-maximal threshold graph with vertex sequence $(a_1,\ldots,a_k;b_1,\ldots,b_k)$,
 $2\le k\le3$. Exclude the case that
$k=2$, $b_1=1$. 
Then
\begin{equation}\label{eq:residual-bound}
 8e(G-B_k)\ge(n-b_k+1)(n-b_k+3).
\end{equation}
\end{lemma}

\begin{proof}
Set $G^*=G-B_k$. There are four cases. 

\noindent {\bf Case 1. $k=2$, $b_1>1$.}
By Lemma \ref{lem:qmax-blocks}(1), $b_1\ge a_1+a_2+1$. Thus 
$
 \nu=a_1+a_2+b_1\le2b_1-1.
$
Since the vertices of $B_1$, together with one vertex of $A_1$, induce
$K_{b_1+1}$, we have
\[
 8e(G^*)\ge4b_1(b_1+1)
       =(2b_1)(2b_1+2)\ge(n-b_k+1)(n-b_k+3).
\]

\noindent {\bf Case 2. $k=3$, $b_1,b_2>1$.}
Applying Lemma \ref{lem:qmax-blocks}(1) at $b_1$ and $b_2$, we have
$
 b_1+b_2\ge a_1+2a_2+a_3+2.
$
Consequently,
\[
n-b_3=a_1+a_2+a_3+b_1+b_2
       \le2(b_1+b_2)-a_2-2\le2(b_1+b_2)-3.
\]
The graph $H$ contains $K_{b_1+b_2+1}$, obtained from
$B_1\cup B_2$ and one vertex of $A_1$. It follows that
\begin{align*}
 8e(G^*)&\ge4(b_1+b_2)(b_1+b_2+1)
    >(2(b_1+b_2)-2)(2(b_1+b_2))
    \ge(n-b_k+1)(n-b_k+3).
\end{align*}

\noindent {\bf Case 3. $k=3$, $b_1=1$.}
Lemma \ref{lem:qmax-singletons}(1) implies $b_2>1$. 
Applying Lemma \ref{lem:qmax-blocks}(2), we have
$
 b_2\ge a_1+a_2+a_3$. Then  $\nu \le 2b_2+1.$
Since $H$ contains $K_{b_2+2}$, formed by $B_2$, the unique
vertex of $B_1$, and one vertex of $A_1$,
\[
 8e(G^*)\ge4(b_2+2)(b_2+1)
       =(2b_2+2)(2b_2+4)\ge(n-b_k+1)(n-b_k+3).
\]

\noindent {\bf Case 4. $k=3$, $b_2=1$.}
Lemma \ref{lem:qmax-singletons}(1) implies  $b_1>1$. 
Applying Lemmas \ref{lem:qmax-blocks}(2) and (3), we have
\begin{equation*}
 b_1\ge a_1+a_2+1,\qquad b_1\ge a_2+a_3+1.
\end{equation*}
The second inequality implies
$n-b_k=a_1+a_2+a_3+b_1+1\le a_1+2b_1.$
The exact residual edge count is
\[
 e(G^*)=\binom{b_1+1}{2}+a_1(b_1+1)+a_2.
\]
Therefore,
$8e(G^*)-(n-b_k+1)(n-b_k+3)
\ge8e(G^*)-(a_1+2b_1+1)(a_1+2b_1+3)=(a_1-1)(4b_1+3-a_1)+8a_2>0.$
The final inequality follows from $a_1,a_2\ge1$ and
$b_1\ge a_1+a_2+1$.

Consequently, the proof is complete. 
\end{proof}

\begin{proof}[Proof of Theorem~\ref{thm:main2}]
Choose a regular $Q$-maximal threshold graph $G$. Applying Theorem \ref{thm:main1}, 
the character of $G$ is at most three. 
One can verify that $G$ is quasi-star and then the result follows if the character  
of $G$ is either
 $k=1$, or $k=2$ and $b_1=1$. Therefore, we assume that the vertex-sequence
  of $G$ is either
 $(a_1,a_2;b_1,b_2)$ with $b_1>1$ or $(a_1,a_2, a_3;b_1,b_2, b_3)$ 
 with $(b_1,b_2)\neq (1,1)$; the latter exclusion is Lemma \ref{lem:qmax-singletons}(1).

In view of the $b_k$ terminal one-vertices are universal, then 
\begin{equation*}
 m=e(G^*)+b_k(n-b_k)+\binom{b_k}{2}.
\end{equation*}
Using \eqref{eq:residual-bound}, we obtain
\begin{align}
 8m
 &\ge(n-b_k+1)(n-b_k+3)+8b_k(n-b_k)+8b_k(b_k-1)\notag\\
 &=n^2+10n-8+(b_k-1)(6n-3b_k-11).
 \label{eq:barrier-factorization}
\end{align}
Here $b_k\ge1$. In every case under consideration the residual graph $G^*$
has at least three vertices, so $b_k\le n-3$. Thus
$
 6n-3b_k-11\ge3n-2>0,
$
and the last product in \eqref{eq:barrier-factorization} is
nonnegative. This proves $8m>n^2+10n-8$.
Then the result follows immediately.
\end{proof}

The bound is sharp.  $\rho(S_{6,11})<\frac{33}{4}<\rho(C_{6,11})=\frac{9+\sqrt{57}}2$ 
and $11=(n^2+10n-8)/8$.

\noindent {\bf Conjecture.}
For every integer pair  $(n,m)$ with
$n\ge 7$ and $n-1\le m\le\frac{n(n-1)}4,$
the quasi-star is the unique   maximal regular threshold graph within $C_{n,m}$.

\noindent {\bf Remark 1.} Let $n=9$, $m=25$, and
\[
 H=C^\triangle_{9,25}=(1,2;\,4,2)=K_2\join(K_5\disj2K_1).
\]
Then
\begin{equation*}
 \max\{\rho(C_{9,25}),\rho(S_{9,25})\}
 <\frac{90}{7}<\rho(H).
\end{equation*}

\noindent {\bf Remark 2.}

\section{Proofs of the results in Section~\ref{sec:qmax-structure}}
\label{sec:proofs-structure}

This section collects the proofs of the results stated
in Section~\ref{sec:qmax-structure}.
 We should emphasize once more that, in the proof of this section,
 zero-blocks or one-blocks of size zero are permitted in temporary block descriptions.
  An internal empty zero-block or one-block can be eliminated by merging the adjacent
   one-blocks or zero-blocks, respectively. 
This operation changes neither the scalar function nor the total number of vertices of either type.

\subsection{Proof of Lemma~\ref{lem:qmax-blocks}}

\begin{proof}
Put $\lambda=\rho(G)\ge n$. For each graph $G'$ constructed below,
let $L'_j(z)$ and $R'_j(z)$ be the functions associated with $G'$,
defined as in \eqref{eq:scalar-LR}. All unlisted block sizes are
left unchanged. We use
\[
S_G(\lambda)=\sum_{j=1}^{k}\frac1{L_j}
\left(\frac1{R_j}-\frac1{R_{j-1}}\right),\qquad
L_i-R_i=\sum_{j<i}(a_j+b_j)+a_i-2.
\]
Moreover, $\lambda\ge n$ gives
$L_i\ge a+\sum_{j<i}b_j$ and $R_i\ge\sum_{j>i}a_j+2$.
These bounds ensure that all denominator factors below are positive
under the corresponding construction conditions.
By Lemma~\ref{lem:scalar-comparison}, it suffices in each case
to obtain $G'\in\Cnm$ with the same totals $a,b$ and
$S_{G'}(\lambda)-S_G(\lambda)>0$.

\noindent\textbf{(i)}
Suppose, to the contrary, that $2\le b_i\le a_i+a_{i+1}$.
Take $r=\min\{a_i,b_i-1\}$, so that
$1\le r\le a_i$ and $1\le b_i-r\le a_{i+1}$.
Delete the $r(b_i-r)$ edges between $r$ chosen vertices of $A_i$
and $b_i-r$ chosen vertices of $B_i$, and add the $r(b_i-r)$ edges
between $b_i-r$ chosen vertices of $A_{i+1}$ and the remaining $r$
vertices of $B_i$. The resulting graph $G'$ has local block data
\[
(a_i,a_{i+1};\,b_i,b_{i+1})
\longmapsto
(a_i-r,b_i,a_{i+1}-b_i+r;\,b_i-r,r,b_{i+1}).
\]
This preserves $n,m,a,b$ and leaves the final one-block unchanged,
so $G'\in\Cnm$. The changed parameters satisfy
\[
\begin{aligned}
L'_i&=L_i, & L'_{i+1}&=L_{i+1}-r, & L'_{i+2}&=L_{i+1},\\
R'_i&=R_i+r, & R'_{i+1}&=R_i-b_i+r, & R'_{i+2}&=R_{i+1}.
\end{aligned}
\]
Cancelling the common terms gives
\[
\begin{aligned}
S_{G'}(\lambda)-S_G(\lambda)
&=-\frac{r}{L_iR_iR'_i}
+\frac{b_i}{L'_{i+1}R'_iR'_{i+1}}
-\frac{b_i-r}{L_{i+1}R_iR'_{i+1}}\\
&=\frac{r}{R_iR'_i}
\left(\frac1{L'_{i+1}}-\frac1{L_i}\right)
+\frac{b_i-r}{R_iR'_{i+1}}
\left(\frac1{L'_{i+1}}-\frac1{L_{i+1}}\right)\\
&=-\frac{r(b_i-r)}{L_iL'_{i+1}R_iR'_i}
+\frac{r(b_i-r)}{L_{i+1}L'_{i+1}R_iR'_{i+1}}
\quad\text{(by }L'_{i+1}-L_i=b_i-r,\;
L_{i+1}-L'_{i+1}=r\text{)}\\
&=\frac{r(b_i-r)}{R_iL'_{i+1}}
\left(\frac1{L_{i+1}R'_{i+1}}-\frac1{L_iR'_i}\right)\\
&=\frac{rb_i(b_i-r)(L'_{i+1}-R_i)}
{L_iL_{i+1}L'_{i+1}R_iR'_iR'_{i+1}}
\quad\text{(by }L_iR'_i-L_{i+1}R'_{i+1}
=b_i(L'_{i+1}-R_i)\text{)}.
\end{aligned}
\]
Now $L'_{i+1}-R_i
=\sum_{j<i}(a_j+b_j)+a_i+b_i-r-2\ge0$;
equality would force $(i,a_i,r,b_i)=(1,1,1,2)$,
contrary to regularity. Hence $L'_{i+1}-R_i>0$ and
$S_{G'}(\lambda)-S_G(\lambda)>0$, proving \textup{(i)}.

\noindent\textbf{(ii)}
Let $b_i=1$. We first exclude $b_{i+1}=1$.
If $b_{i+1}=1$, choose one vertex from $A_i$
and one from $A_{i+2}$. Delete the edge between the vertex of $A_i$
 and the vertex of $B_i$, and add the edge between the vertex of $A_{i+2}$
 and the vertex of $B_{i+1}$. This gives $G'\in\Cnm$
with local replacement
\[
(a_i,a_{i+1},a_{i+2};\,1,1,b_{i+2})
\longmapsto
(a_i-1,a_{i+1}+2,a_{i+2}-1;\,1,1,b_{i+2}).
\]
Here all $L'_j=L_j$, while the only changed $R$-parameters are
$R'_i=R_i+1$ and $R'_{i+1}=R_{i+1}-1$.
Hence,
\[
\begin{aligned}
S_{G'}(\lambda)-S_G(\lambda)
&=\left(\frac1{L_i}-\frac1{L_{i+1}}\right)
\left(\frac1{R'_i}-\frac1{R_i}\right)
+\left(\frac1{L_{i+1}}-\frac1{L_{i+2}}\right)
\left(\frac1{R'_{i+1}}-\frac1{R_{i+1}}\right)\\
&=\frac1{L_iL_{i+1}}\left(\frac1{R'_i}-\frac1{R_i}\right)
+\frac1{L_{i+1}L_{i+2}}
\left(\frac1{R'_{i+1}}-\frac1{R_{i+1}}\right)
\;\text{(by }L_{i+1}-L_i=L_{i+2}-L_{i+1}=1\text{)}\\
&=-\frac1{L_iL_{i+1}R_iR'_i}
+\frac1{L_{i+1}L_{i+2}R_{i+1}R'_{i+1}}
\quad\text{(by }R'_i-R_i=R_{i+1}-R'_{i+1}=1\text{)}.
\end{aligned}
\]
Now $a_{i+1}\ge1$, $L_i\ge R_i-1$, and $R_i>2$ imply
$L_{i+2}R_{i+1}R'_{i+1}
\le L_i(R_i+1)(R_i-2)<L_iR_iR'_i$.
Hence $S_{G'}(\lambda)-S_G(\lambda)>0$, a contradiction.
Thus $b_{i+1}>1$.

Suppose now that $b_{i+1}+1\le a_i+a_{i+1}+a_{i+2}$.
By \textup{(i)},
$a_{i+1}+a_{i+2}<b_{i+1}<a_i+a_{i+1}+a_{i+2}$.
Make the local replacement
\begingroup
\medmuskip=1.5mu
\thickmuskip=3mu
\[
(a_i,a_{i+1},a_{i+2};1,b_{i+1},b_{i+2})
\longmapsto(a_i+a_{i+1}+a_{i+2}-b_{i+1}-1,b_{i+1},1;
a_{i+2},a_{i+1},b_{i+1}+1-a_{i+1}-a_{i+2}+b_{i+2}).
\]
\endgroup
With the zero-vertices and one-vertices matched in their original
orders, this exchange deletes and adds
$a_{i+1}(a_{i+2}-1)+a_{i+2}(b_{i+1}-a_{i+1}-a_{i+2}+1)$
edges in each direction. All new block sizes are nonnegative,
the final one-block remains positive, and $n,a,b$ are preserved.
Thus the resulting graph $G'$ belongs to $\Cnm$.
Its changed parameters satisfy
\[
\begin{aligned}
L'_i&=L_i, &
L'_{i+1}&=L_i+a_{i+2}, &
L'_{i+2}&=L_i+a_{i+1}+a_{i+2},\\
R'_i&=R_i+b_{i+1}-a_{i+1}-a_{i+2}+1, &
R'_{i+1}&=R_{i+2}+1, &
R'_{i+2}&=R_{i+2}.
\end{aligned}
\]
Cancelling the common terms and regrouping gives
\[
\begin{aligned}
S_{G'}(\lambda)-S_G(\lambda)
&=\frac{a_{i+2}}{L_iL'_{i+1}R'_i}
+\frac{a_{i+1}}{L'_{i+1}L'_{i+2}R'_{i+1}}
+\frac{R'_i-R_i}{L'_{i+2}L_{i+2}R_{i+2}}
-\frac1{L_iL_{i+1}R_i}
-\frac{b_{i+1}}{L_{i+1}L_{i+2}R_{i+1}}\\
&=\left[
\frac{(a_{i+2}-1)(b_{i+1}-a_{i+2}+1)}{L'_{i+1}R_{i+1}}
+\frac{R'_i-R_i}{L'_{i+2}R_{i+2}}
\right]
\left(\frac1{L_{i+2}R'_{i+1}}-\frac1{L_{i+1}R'_i}\right)\\
&\quad+\frac{R'_i-R_i}{L_{i+1}R'_i}
\left(\frac1{L'_{i+2}R_{i+2}}-\frac1{L_iR_i}\right).
\end{aligned}
\]
Now $L_i-R_{i+2}
=\sum_{j<i}(a_j+b_j)+a_i+a_{i+1}+a_{i+2}-2>0$, so
$L_{i+1}R'_i-L_{i+2}R'_{i+1}
=b_{i+1}(L_i-R_{i+2})>0$ and 
$L_iR_i-L'_{i+2}R_{i+2}
=(a_{i+1}+a_{i+2})(L_i-R_{i+2})>0.$
Together with 
$R'_i-R_i=b_{i+1}-a_{i+1}-a_{i+2}+1>0$,
$a_{i+2}\ge1$, and $b_{i+1}>a_{i+1}+a_{i+2}$,
we obtain $S_{G'}(\lambda)-S_G(\lambda)>0$.
This proves \textup{(ii)}.

\noindent\textbf{(iii)}
Let $b_i>1$ and $b_{i+1}=1$, and suppose, to the contrary, that
$b_i\le a_{i+1}+a_{i+2}$.
By \textup{(i)}, $b_i>a_i+a_{i+1}>a_{i+1}$.
Choose one vertex from $A_i$, $b_i-a_{i+1}$ vertices from $B_i$,
and $b_i-a_{i+1}$ vertices from $A_{i+2}$.
Delete the edges between the chosen vertex of $A_i$ and the chosen
vertices of $B_i$, and add the edges between the chosen vertices
of $A_{i+2}$ and the vertex of $B_{i+1}$.
The resulting graph $G'\in\Cnm$ has local block data
\[
(a_i,a_{i+1},a_{i+2};\,b_i,1,b_{i+2})
\longmapsto(a_i-1,1,b_i,a_{i+1}+a_{i+2}-b_i;\,
b_i-a_{i+1},a_{i+1},1,b_{i+2}).
\]
The changed parameters satisfy
\[
\begin{aligned}
L'_i&=L_i, &
L'_{i+1}&=L_i+b_i-a_{i+1}, &
L'_{i+2}&=L_{i+1}, &
L'_{i+3}&=L_{i+2},\\
R'_i&=R_i+1, &
R'_{i+1}&=R_i, &
R'_{i+2}&=R_i-b_i, &
R'_{i+3}&=R_{i+2}.
\end{aligned}
\]
Hence
\[
\begin{aligned}
S_{G'}(\lambda)-S_G(\lambda)
&=\frac1{L_{i+1}L_{i+2}}
\left(\frac1{R'_{i+2}}-\frac1{R_{i+1}}\right)
+\frac{b_i-a_{i+1}}{L_iL'_{i+1}}
\left(\frac1{R'_i}-\frac1{R_i}\right)\\
&=(b_i-a_{i+1})
\left[
\frac1{L_{i+1}L_{i+2}R'_{i+2}R_{i+1}}
-\frac1{L_iL'_{i+1}R_iR'_i}
\right]
\quad\text{(by }R_{i+1}-R'_{i+2}=b_i-a_{i+1}\text{)}.
\end{aligned}
\]
Now $L_{i+1}-R_i=\sum_{j<i}(a_j+b_j)+a_i+b_i-2>0$, so
$L_iR_i-L_{i+1}R'_{i+2}=b_i(L_{i+1}-R_i)>0$ and $L'_{i+1}R'_i-L_{i+2}R_{i+1}
=(a_{i+1}+1)(L_{i+1}-R_i)>0.$
Since $b_i-a_{i+1}>0$, we obtain $S_{G'}(\lambda)-S_G(\lambda)>0$.
This proves \textup{(iii)} and completes the proof.
\end{proof}

\subsection{Proof of Lemma~\ref{lem:qmax-singletons}}

\begin{proof}
(i). Suppose first that $b_i=b_j=1$ with
$j\ge i+2$ and $j<k$.  Let $G_+$ and $G_-$ be obtained by replacing $
 (a_i,a_{i+1},a_j,a_{j+1})$ by $
 (a_i+1,a_{i+1}-1,a_j-1,a_{j+1}+1)$
and $
 (a_i-1,a_{i+1}+1,a_j+1,a_{j+1}-1),$
respectively.  Any empty blocks are suppressed.  Both $G_+$ and $G_-$ preserve
$n$ and $m$.  From \eqref{eq:positional-S},
\begin{align*}
 S_{G_+}(\lambda)-S_G(\lambda)
 &=\frac{\omega_{s_i+1}}{L_i(L_i+1)}
   -\frac{\omega_{s_j}}{L_j(L_j+1)},\\
 S_{G_-}(\lambda)-S_G(\lambda)
 &=-\frac{\omega_{s_i}}{L_i(L_i+1)}
   +\frac{\omega_{s_j+1}}{L_j(L_j+1)}.
\end{align*}
Then  $S_{G_+}(\lambda)+S_{G_-}(\lambda)-2S_G(\lambda)>0$ because
$\omega_\ell(\lambda)$ is strictly increasing on $\ell$. 
 Hence $\max \{\rho(G_+),\rho(G_+)\}>\rho(G)$, contradicting
contradicting the $Q$-maximality of $G$
by Lemma~\ref{lem:scalar-comparison}.  The adjacent case is excluded by the first exchange in the
proof of Lemma~\ref{lem:qmax-blocks}\textup{(ii)},
which does not use regularity.

(ii). Suppose that $a_i=b_j=1$ for indices in the stated
range.  Form $G_+$ and $G_-$ by changing $
  (b_{i-1},b_i;a_j,a_{j+1})$
to $
  (b_{i-1}-1,b_i+1;a_j-1,a_{j+1}+1)$
and $ (b_{i-1}+1,b_i-1;a_j+1,a_{j+1}-1),$
respectively.  Moving one vertex from $B_{i-1}$ to $B_i$ changes the
edge count by $a_i=1$, while moving one vertex from $A_j$ to
$A_{j+1}$ changes it by $-b_j=-1$.  Since
$j\notin\{i-1,i\}$, there is no cross term.  Thus both of $G_+$ and $G_-$
belong to $\mathcal C_{n,m}$.
At $\lambda=\rho(G)$ the two scalar differences are
\begin{align*}
 S_{G_+}(\lambda)-S_G(\lambda)
 &=\omega_{s_i}\left(\frac1{L_i-1}-\frac1{L_i}\right)
   +\omega_{s_j}\left(\frac1{L_j+1}-\frac1{L_j}\right),\\
 S_{G_-}(\lambda)-S_G(\lambda)
 &=\omega_{s_i}\left(\frac1{L_i+1}-\frac1{L_i}\right)
   +\omega_{s_j+1}\left(\frac1{L_j}-\frac1{L_j+1}\right).
\end{align*}
Since $\lambda\ge n$ gives $R_j\ge2$, and  $L_i\ge2$, then
\begin{equation*}S_{G_+}(\lambda)+S_{G_-}(\lambda)-2S_G(\lambda)=\frac{2}{R_i(R_i+1)L_i(L_i^2-1)}
 +\frac{2}{R_j(R_j^2-1)L_j(L_j+1)}>0.
\end{equation*}
Hence, $\max \{\rho(G_+),\rho(G_+)\}>\rho(G)$, again a
contradiction.
\end{proof}

\subsection{Proof of Lemma~\ref{lem:strict-midpoint} }

\begin{proof}
Let $G_+$ and $G_-$ be  graphs obtained from $G$ by replacing
$b_i$ with $b_i+\delta_i$ and $b_i-\delta_i$, respectively.  The two
displayed moment conditions preserve the number of one-vertices and the
number of edges, so $G_+,G_-\in\mathcal C_{n,m}$.

Put $z=\rho(G)$ and $
 d_i=\sum_{j\ge i}\delta_j.$
The functions $R_i(z)$ are unchanged, while the two denominators that
replace $L_i(z)$ are $L_i(z)-d_i$ and $L_i(z)+d_i$.  Hence
\begin{align*}
 S_{G_+}(z)+S_{G_-}(z)-2S_G(z)
 &=
 \sum_{i=1}^k
 \frac{a_i}{R_{i-1}(z)R_i(z)}
 \left(
  \frac1{L_i(z)-d_i}+
  \frac1{L_i(z)+d_i}-
  \frac2{L_i(z)}
 \right)\\
 &=
 \sum_{i=1}^k
 \frac{2a_i d_i^2}
 {R_{i-1}(z)R_i(z)L_i(z)\bigl(L_i(z)^2-d_i^2\bigr)}>0.
\end{align*}
All denominators are positive because both perturbed words are valid.
Since $(\delta_i)$ is nonzero and has zero total, at least one tail sum
$d_i$ is nonzero.  Therefore one of $G_+$ and $G_-$ satisfies
$S_{G_\pm}(z)>S_G(z)$, and Lemma~\ref{lem:scalar-comparison} gives
$\max \{\rho(G_+),\rho(G_-)\}>\rho(G)$, as proved.
\end{proof}

\subsection{Proof of Lemma~\ref{lem:splitting-inequality} }

\begin{proof}
For $z\ge n$, by  a direct substitution, we have
\begin{align*}
 S_{G'}(z)-S_G(z)
 &=\frac{r}
 {[L_i(z)-\tau(a_i-r)][R_i(z)+a_i][R_i(z)+a_i-r]}\\
 &\quad+\frac{a_i-r}
 {[L_i(z)+\tau r][R_i(z)+a_i-r]R_i(z)}-\frac{a_i}
 {L_i(z)R_i(z)[R_i(z)+a_i]}\\
 &=\frac{\tau a_i r(a_i-r)
 \{\tau[R_i(z)+a_i-r]-L_i(z)\}}{\Phi(z)},
\end{align*}
where
\begin{align*}
 \Phi(z)={}&L_i(z)R_i(z)[R_i(z)+a_i][R_i(z)+a_i-r][L_i(z)-\tau(a_i-r)][L_i(z)+\tau r]>0.
\end{align*}
By the maximality of $G$ and Lemma
\ref{lem:scalar-comparison}, $
  S_{G'}(\rho(G))-S_G(\rho(G))\le0.$
Since $\tau a_i r(a_i-r)>0$, inequality
\eqref{eq:splitting-bound} follows.  Taking $\tau=2$ and
$r=\lfloor a_i/2\rfloor$ gives
\eqref{eq:splitting-special-bound}.
\end{proof}

\subsection{Proof of Lemma~\ref{lem:splitting-consequences}}

\begin{proof}
Let $H$ be the threshold graph with vertex-sequence $
 (a_1,\ldots,a_{i-1};\,b_1,\ldots,b_{i-1}),
 $ and put $
 t_i=\sum_{j=i}^k b_j.$
Since $H\join K_{t_i}$ is a subgraph of $G$, Lemmas
\ref{lem:clique-join-estimate} and \ref{lem:prefix-estimates} give
\begin{equation}\label{eq:split-lower-rho}
 \rho(G)>
 \frac32\sum_{j<i}(a_j+b_j)+2t_i-2.
\end{equation}
Suppose that an integer $r$ satisfies \eqref{eq:split-range} with
$\tau=3$.  Then \eqref{eq:splitting-bound} gives
\begin{equation}\label{eq:split-upper-rho}
 2\rho(G)\le
 3\sum_{j<i}(a_j+b_j)+2\sum_{j=i}^k b_j+3r-6.
\end{equation}
Combining \eqref{eq:split-lower-rho} and
\eqref{eq:split-upper-rho} yields $
 0>2\sum_{j=i}^k b_j-3r+2
 \ge\sum_{j=i}^k b_j+2>0,$
a contradiction.  Therefore the interval in
\eqref{eq:split-range} is empty.

By Lemma~\ref{lem:qmax-blocks}\textup{(i)},
$b_{i-1},b_i\ge a_i+2\ge4$, so both floor values below are at least
one.  Since $a_i\ge2$, emptiness of the interval is equivalent to $
 a_i-\left\lfloor\frac{b_{i-1}}3\right\rfloor
 >\left\lfloor\frac{b_i}3\right\rfloor,$
and hence
\[
 \left\lfloor\frac{b_{i-1}}3\right\rfloor
 +\left\lfloor\frac{b_i}3\right\rfloor\le a_i-1.
\]
It follows that
\begin{align*}
 b_{i-1}+b_i\le
 \left(3\left\lfloor\frac{b_{i-1}}3\right\rfloor+2\right)
 +\left(3\left\lfloor\frac{b_i}3\right\rfloor+2\right)\le3(a_i-1)+4=3a_i+1.
\end{align*}
Finally, Lemma~\ref{lem:qmax-blocks}\textup{(i)} gives $
 3a_i+1\ge b_{i-1}+b_i
 \ge2a_i+a_{i-1}+a_{i+1}+2,$
so $a_i\ge a_{i-1}+a_{i+1}+1$.
\end{proof}

\subsection{Proof of Lemma~\ref{lem:clique-join-estimate} }

\begin{proof}
Let $\mathbf x$ be a nonnegative unit Perron vector of $Q(G)$.  The complement identity gives
\[
  (\one^{\mathsf T}\mathbf x)^2
  =\rho(G)-n+2+
    \mathbf x^{\mathsf T}Q(\overline G)\mathbf x
  \ge \rho(G)-n+2.
\]
Set $M=(\rho(G)+t)I-Q(G).$
Since every eigenvalue of $M$ is at least $t>0$ and $M\mathbf x=t\mathbf x$,
the matrix $M$ is positive definite.
Applying the Cauchy--Schwarz inequality, 
\begin{align*}
 (\one^{\mathsf T}\mathbf x)^2
 &=\bigl((M^{-1/2}\one)^{\mathsf T}M^{1/2}\mathbf x\bigr)^2\le(\one^{\mathsf T}M^{-1}\one)
       (\mathbf x^{\mathsf T}M\mathbf x)
  =t\Gamma_G(\rho(G)+t).
\end{align*}
For $z>\rho(G)+t$, put
\[
  f_{t,G}(z)=z-n-2t+2-t\Gamma_G(z-t).
\]
The Schur-complement equation for $Q(G\join K_t)$ shows that its
largest zero is $\rho(G\join K_t)$.  At $z=\rho(G)+2t$,
\[
 f_{t,G}(\rho(G)+2t)
 =\rho(G)-n+2-t\Gamma_G(\rho(G)+t)\le0.
\]
Hence $\rho(G\join K_t)\ge\rho(G)+2t$.
\end{proof}

\subsection{Proof of Lemma~\ref{lem:rayleigh-lower-bounds} }

\begin{proof}
By the Rayleigh principle, it suffices to find a nonzero vector
$\mathbf v$ such that
$\mathbf v^{\mathsf T}Q(G)\mathbf v
  -(\frac32n-2)\|\mathbf v\|^2>0.$

\noindent\textbf{(1)}
Assign coordinates $1/2$ and $1$ to the zero- and
one-vertices, respectively.  Then
\begin{align*}
 \mathbf v^{\mathsf T}Q(G)\mathbf v
 -\left(\frac32(a+b)-2\right)\|\mathbf v\|^2 &~\ge
 \frac94a^2+2b(b-1)
 -\left(\frac32(a+b)-2\right)\left(\frac a4+b\right)\\
 &~=
 \frac18\left\{
 4\left(b-\frac{15a}{8}\right)^2
 +\frac{15a^2}{16}+4a
 \right\}>0.
\end{align*}

\noindent\textbf{(2)} Assign coordinates $2,5,1,8$ to the zero- and
one-vertices of the first and second subsequences, respectively.
Then $
  \|\mathbf v\|^2=4a+25b+c+64d$
and thus
\[
 \mathbf v^{\mathsf T}Q(G)\mathbf v
 \ge49a^2+50b(b-1)+81c^2+128d(d-1)+100ad+169bd.
\]
Consequently,
\begin{align*}
 &~~~~~~~\mathbf v^{\mathsf T}Q(G)\mathbf v
 -\left(\frac32(a+b+c+d)-2\right)\|\mathbf v\|^2\\
 &\quad\ge
 \frac{(56c-26a-7(b-a)-67(d-c))^2}{224}
 +\frac{161(8a-9(b-a))^2}{1152}\\
 &\qquad
 +\frac{19a^2}{504}
 +\frac{123(b-a)^2}{128}
 +\frac{2679(d-c)^2}{224}
 +\frac{1005a(d-c)}{56}
 +\frac{501(b-a)(d-c)}{16}+8a+2c>0.
\end{align*}
\end{proof}

\subsection{Proof of Lemma~\ref{lem:prefix-estimates} }
\begin{proof}
If $b_1>1$, then, from Lemma~\ref{lem:qmax-blocks}\textup{(i)}, we have $
 a_1b_1\ge a_1(a_1+a_2+1)>a_1^2,$
so \textup{(1)} follows from Lemma
\ref{lem:rayleigh-lower-bounds}\textup{(1)}.

For \textup{(2)}, suppose first that $b_1>1$.  By Lemma
\ref{lem:qmax-blocks}\textup{(i)},
\[
 a_1b_1+(a_1+a_2)b_2
 \ge a_1(a_1+a_2)+(a_1+a_2)a_2
 =(a_1+a_2)^2.
\]
If $b_1=1$, Lemma~\ref{lem:qmax-blocks}\textup{(ii)} implies
\[
 a_1b_1+(a_1+a_2)b_2
 \ge a_1+(a_1+a_2)^2>(a_1+a_2)^2.
\]
In both cases the conclusion follows from Lemma
\ref{lem:rayleigh-lower-bounds}\textup{(1)}.
\end{proof}

\subsection{Proof of Lemma~\ref{lem:revised-prefix-bounds} }
\begin{proof}
\noindent\textup{(1)}
Let $\mathbf x$ be constant on $A_1,B_1,A_2,B_2$
with respective values $(2,3,0,4)$.
The assumptions imply $a_1\le a_2+1$ and $b_1\le2a_2$.
A direct calculation gives
\[
b_1+4b_2-4a_1-2a_2+8\ge12,~~
4b_2-a_1-3a_2+6\ge13,~~
4b_2-b_1-2a_1+8\ge14.
\]
These are the respective lower bounds for the coordinates
of $[Q(F)-(2b_1+b_2+2a_1+a_2-4)I]\mathbf x$
on $A_1,B_1,B_2$.
Since $\mathbf x$ vanishes on $A_2$, we obtain
\[
\mathbf x^{\mathsf T}Q(F)\mathbf x-(2b_1+b_2+2a_1+a_2-4)\|\mathbf x\|^2\ge24a_1+39b_1+56b_2>0.
\]
The Rayleigh principle proves \textup{(1)}.

\smallskip
\noindent\textup{(2)}
Let $F_0$ have
vertex-sequence
$(a_1,1,a_3;\,b_1,b_2,a_3+2).$
Then $F=F_0\vee K_t$, and
Lemma~\ref{lem:clique-join-estimate} gives
$\rho(F)\ge\rho(F_0)+2t,$ where  $t=b_3-a_3-2\ge0$,
with the case $t=0$ being immediate.
Since the proposed lower bound increases by only $t$,
it suffices to consider $b_3=a_3+2$.

Let $\mathbf x$ be constant on
$A_1,B_1,A_2,B_2,A_3,B_3$
with respective values $(5,7,0,7,1,9)$.
The assumptions imply
$a_1\ge1$, $b_1\ge3$, $a_3\ge4$, and $b_2\ge6$.
Thus,
\[
(b_1-a_1-2)(50a_1+54b_1+18a_3+4b_2-193)\ge0,
~~
(a_3-b_1-1)(54b_1+78a_3+4b_2-494)\ge0.
\]
Expanding the Rayleigh difference gives the sum of these
two nonnegative terms and
\[
10a_3(2a_3-b_2)+200(b_2-a_3-2)+1333a_3+168.
\]
Since $a_3+2\le b_2\le2a_3$, it follows that
\[
\mathbf x^{\mathsf T}Q(F)\mathbf x-(2b_1+2b_2+b_3+2a_1+a_3-3)\|\mathbf x\|^2\ge1333a_3+168>0.
\]
The Rayleigh principle proves \textup{(2)}.
\end{proof}

\section*{Acknowledgments}

This work was supported by the National Natural Science Foundation of China (No.  12271484).

%%%%%%%%%%%%%%%%%%%%%%%%%%%%
\section*{Declaration of interests}
%%%%%%%%%%%%%%%%%%%%%%%%%%%%

The authors declare that they have no conflict of interest.

%%%%%%%%%%%%%%%%%%%%%%%%%%%%
\section*{Data availability}
%%%%%%%%%%%%%%%%%%%%%%%%%%%%

Our manuscript has no associated data.

\section*{Declaration on the use of generative AI}
During the preparation of this manuscript, the authors used the AI tool ChatGPT
 to assist with language polishing and organizing literature identified by the authors.
  The authors independently verified all references, reviewed and revised the AI-assisted
   content, and are responsible for the accuracy,
 originality, and integrity of the final manuscript.

\end{document}